\documentclass[twoside,11pt]{article}

\usepackage[round, authoryear]{natbib} 

\usepackage{jmlr2e}

\usepackage{amsmath}
\usepackage[mathscr]{euscript}
\usepackage{bbold}

\usepackage{graphicx}
\graphicspath{ {images/} }
\usepackage{caption}
\usepackage{subcaption}

\usepackage{booktabs}
\usepackage[table,xcdraw]{xcolor}
\usepackage[normalem]{ulem}
\useunder{\uline}{\ul}{}

\ShortHeadings{Non-Parametric Test for IV}{Helder R., Cirilo A. and Nilton R.}
\firstpageno{1}

\begin{document}

\title{Statistical Hypothesis Testing  for Information Value (IV)}

\author{\name Helder Rojas
\email  helder\_rojas@hotmail.com \\
       \addr  Imperial College London \\
       South Kensington Campus \\
       London, England
       \AND
       \name Cirilo Alvarez \email calvarezr@uni.edu.pe \\
       \addr  Professional School of Statistical Engineering (UNI)\\
       National Engineering University\\
       Lima, Perú
       \AND
       \name Nilton Rojas \email nrojasv@uni.pe \\
       \addr  Professional School of Computer Science\\
       National Engineering University (UNI)\\
       Lima, Perú
       }

\editor{}

\maketitle

\begin{abstract}
Information Value (IV) is a widely used technique for feature selection prior to the modeling phase, particularly in credit scoring and related domains. However, conventional IV-based practices rely on fixed empirical thresholds, which lack statistical justification and may be sensitive to characteristics such as class imbalance. In this work, we develop a formal statistical framework for IV by establishing its connection with Jeffreys divergence and propose a novel nonparametric hypothesis test, referred to as the J-Divergence test. Our method provides rigorous asymptotic guarantees and enables interpretable decisions based on \(p\)-values. Numerical experiments, including synthetic and real-world data, demonstrate that the proposed test is more reliable than traditional IV thresholding, particularly under strong imbalance. The test is model-agnostic, computationally efficient, and well-suited for the pre-modeling phase in high-dimensional or imbalanced settings. An open-source Python library is provided for reproducibility and practical adoption.
\end{abstract}
\hfill \break

\section{Introduction}\label{intro}
Feature selection in classification problems, or also known as predictor variable reduction, is an important step before venturing to build a model. In high-dimensional data, data sets with many features, this is even more important, because parsimonious and easy-to-interpret models are prioritized. Information Value (IV) is a quite popular tool for features selection in binary classification problems. Initially, the IV, in conjunction with Weight of Evidence (WoE), inspired by Information Theory developed in the 1940s \citet{howard1966information, kullback1997information}, was designed to selection, raking and monitor predictor variables in the credit risk industry \citet{siddiqi2012credit}. However, in recent years it has received increasing attention in diverse areas such as; marketing, engineering and geosciences \citet{matthijsen2015information, batar2021landslide, alsabhan2022landslide}. In machine learning, in a very useful and interesting way, IV is used to the optimal discretization of a variable into bins given a discrete or continuous numeric target, process known as optimal binning \citet{navas2020optimal}. Furthermore, the implementation using the IV is simple and practical \citet{lund2012collapsing, lin2013variable}.

In order to illustrate how the IV is calculated, consider a binary target $y$ and $x$ a feature for which we want to evaluate its predictive power. The predictor variable $x$ is organized in $r$ bins (groups or buckets). For binning techniques in detail, see \citet{navas2020optimal}. Let $g_j$ the count of instances with $y=1$ at bin $j$, and $g_t$ the count of total instances with $y=1$. Define the empirical relative frequencies $\hat{p}_j$ for each bin $j$ as $\hat{p}_j=\frac{g_j}{g_t}$. Correspondingly, define the empirical relative frequencies $\hat{q}_j=\frac{b_j}{b_t}$ for each bin $j$, where $b_j$ the count of instances with $y=0$ at bin $j$, and $b_t$ the count of total instances with $y=0$. Then, the Weight of Evidence (WoE) for bin $j$ is defined as
begin{document}

\begin{table}[]
\centering 
\begin{tabular}{cc} 
\toprule
IV               & Predictive Power                   \\ 
\midrule
\rowcolor[HTML]{38FFF8} 
$<$ 0.02   & Not useful for prediction      \\
\rowcolor[HTML]{DAE8FC} 
0.02 to 0.1      & Weak                     \\
\rowcolor[HTML]{38FFF8} 
0.1 to 0.3       & Medium                   \\
\rowcolor[HTML]{DAE8FC} 
0.3 to 0.5       & Strong                   \\
\rowcolor[HTML]{38FFF8} 
$>$ 0.5 & Suspicious Predictive Power \\ 
\bottomrule
\end{tabular}
\caption{Threshold values for IV and their corresponding predictive power categories, see \cite{siddiqi2012credit}.} 
\label{Theshold_IV}
\end{table}
\begin{equation*}
    \textrm{WoE}_j=\ln\Big(\frac{\hat{p}_j}{\hat{q}_j}\Big), 
\end{equation*}
and IV of predictor $x$ is defined by 
\begin{equation}\label{IV}
    \textrm{IV}=\sum\limits_{j=1}^{r}(\hat{p}_j-\hat{p}_j)\ln\Big(\frac{\hat{p}_j}{\hat{q}_j}\Big)=\sum\limits_{j=1}^{r}(\hat{p}_j-\hat{p}_j)\cdot\textrm{WoE}_j.
\end{equation}
From Equation \eqref{IV}, informally we can think that $\hat{p}_j-\hat{q}_j$ is the deviation between $\hat{p}_j$ and $\hat{q}_j$ and $\textrm{WoE}_j$ is the importance or weight this deviation. Therefore, we can informally interpret IV as the weighted deviation between the empirical distributions $(\hat{p}_j)_{1\leq j\leq r}$ and $(\hat{q}_j)_{1\leq j\leq r}$. Consequently, IV is used to measure the predictive power, or also called discrimination power, of the a feature in binary classification problems $x$. 

\subsection{Our critique and motivation}
Curiously, there are practical rules, quite familiar but mysterious at the same time, to decide if the IV of a predictor feature $x$ is high enough to be taken into account in the modeling phase, see Table \ref{Theshold_IV}. A fairly popular practice among IV users is to filter variables whose predictive power is medium or strong, that is, variables whose IV is greater than or equal to $0.1$. This selection criterion is simple and eventually leads to good results. However, the thresholds established for the IV seem to be magic and immutable numbers, fixed thresholds, regardless of the possible particularities that may exist in the data. For example, if there is a considerable imbalance in the data with respect to the target variable, a very recurring context in classification problems, the thresholds should be established based on the level of imbalance, and not be the same for cases where there is no such imbalance. Motivated by our intuitive impression about a close connection between the particularity of the data and the optimal thresholds for the IV, we made an exhaustive search about the origins and arguments about the fixed thresholds determined in Table \ref{Theshold_IV}. To our surprise, we have not found anything, any mathematical or theoretical construction, that supports and guarantees the validity this fix thresholds. As we said before, under appropriate conditions, this thresholds work relatively well, we think that these were determined based on the accumulation of empirical experiences. This absence of convincing theoretical arguments motivated us to work on establishing a coherent and precise statistical connection between the predictive power of a variable and existing particularities in the target variable. With this objective in mind, since \eqref{IV} is a sample quantity, it is necessary to define IV as a statistical estimator and study its asymptotic properties and from which it is possible that we can establish a statistical hypothesis test.
\subsection{Main contributions}
In this work we propose a non-parametric hypothesis test for the IV. Taking advantage of its relationship with the divergence measures developed around Information Theory, we study the sampling properties of the IV, statistical estimator properties, in terms of the Jeffreys divergence measure. In Section \ref{test}, inspired in works \citet{basharin1959statistical, stewart2019jensen, ba2019divergence}, we prove Strong Law of Large Numbers (almost sure consistency) and Central Limit Theorem (asymptotic normality) for \eqref{IV}, with which we derived a statistical hypothesis test for this estimator. This test is non-parametric since it does not make any specific assumptions about the shape of the empirical distributions evaluated. In Section \ref{simulation}, through a simulation study we show the efficiency and superiority of our test, we compare its performance with the widely used criterion $\textrm{IV}>0.1$. Additionally, in Section \ref{aplication}, we illustrate the use of our test in fraud data. Finally, in Section \ref{python}, we present a Python library where we have implemented all our results.

\section{Definition of the hypotheses test}\label{test}
Let $S= \{a_1, a_2, \cdots ,a_r\}$ the bins set and $(X,Y) \in S \times\{0,1\}$ be a bivariate random vector endowed with an unknown distribution $\mathbb{P}_{X,Y}=\mathbb{P}$. The random variable $X$ is the predictor feature and $Y$ is the binary target variable. Let $\mathbf{p}=(p_j)_{1\leq j \leq r}$ and $\mathbf{q}=(q_j)_{1\leq j \leq r}$ two discrete probability distributions on $S$, defined as
\begin{equation}\label{p_q}
    p_j=\mathbb{P}(X \in a_j|Y=1)\quad\textrm{and}\quad q_j=\mathbb{P}(X \in a_j|Y=0),
\end{equation}
where $j\in D=\{1,2,\cdots,r\}$. Furthermore, we assume that 
\begin{equation}\label{assumption}
    p_j>0\quad\textrm{and}\quad q_j>0\quad\textrm{for all}\quad j\in D.
\end{equation}
\noindent
In order to measure the discrepancy between the two probability distributions $\mathbf{p}$ and $\mathbf{q}$, we introduce the Jeffreys divergence.
\\

\noindent
{\bf Definition 1} {\it The Jeffreys divergence between the two probability distribution $\mathbf{p}$ and $\mathbf{q}$ is given by
\begin{equation}\label{J_divergence}
    J(\mathbf{p},\mathbf{q})=\sum_{j\in D}(p_j-q_j)\ln\frac{p_j}{q_j}.
\end{equation}
}
\begin{remark}
    The Jeffreys divergence of $\mathbf{p}$ and $\mathbf{q}$ is nonnegative, and equal to 0 if and only if $\mathbf{p}=\mathbf{q}$. Furthermore, $J$ is symmetric, that is, $J(\mathbf{p}, \mathbf{q})=J(\mathbf{q}, \mathbf{p})$. For more details about this measure of divergence, see \citet{jeffreys1946invariant}, \citet{kullback1997information} and \citet{jeffreys1998theory}.
\end{remark}

Since the goal is to evaluate the predictive power of $X$ on $Y$, then we are interested in testing the hypothesis
\begin{equation}\label{HT}
    H_o: \mathbf{p}=\mathbf{q}\quad\textrm{versus}\quad H_A:\mathbf{p}\neq\mathbf{q}.
\end{equation}
Given a sample of $N$ i.i.d. random vectors $(X_1,Y_1) \cdots (X_N,Y_N)$ from the probability distribution $\mathbb{P}$, we define the empirical probability distributions 
\begin{equation}\label{empirical_p}
    {\hat{\mathbf{p}}}_n=\big({\hat{p}}_n^{\,a_j}\big)_{a_j \in S}\quad\textrm{where}\quad\hat{p}_{a_j, n} = \frac{1}{n} \sum_{i=1}^{N} \mathbf{I}[X_i \in a_j] \, \mathbf{I}[Y_i = 1],
\end{equation}
and
\begin{equation}\label{empirical_q}
    {\hat{\mathbf{q}}}_m=\big({\hat{q}}_m^{\,a_j}\big)_{a_j \in S}\quad\textrm{where}\quad\hat{p}_{a_j, n} = \frac{1}{n} \sum_{i=1}^{N} \mathbf{I}[X_i \in a_j] \, \mathbf{I}[Y_i = 0].
\end{equation}
Here $n=\sum\limits_{j=1}^{r}\sum_{i=1}^{N} \mathbf{I}[X_i \in a_j] \, \mathbf{I}[Y_i = 1]$ (instances with label $1$) and $m=N-n$ (instances with label $0$). Using Equations \eqref{empirical_p} and \eqref{empirical_q}, we can directly estimate the Jeffreys divergence \eqref{J_divergence}.
\\

\noindent
{\bf Definition 2}\label{Jeffreys_estimator} {\it The Jeffreys divergence estimator is defind as
\begin{equation}\label{divergence_empirical}
    J(\hat{\mathbf{p}},\hat{\mathbf{q}})=\sum_{j\in D}({\hat{p}}_n^{\,a_j}-{\hat{q}}_m^{\,a_j})\ln\frac{{\hat{p}}_n^{\,a_j}}{{\hat{q}}_m^{\,a_j}}.
\end{equation}
}
\begin{remark}
    Note that IV defined in Equation \eqref{IV} is equivalent to the Jeffreys divergence estimator this definition. Henceforth, we will use this notation to refer to IV.
\end{remark}

The empirical divergence in Equation  \eqref{divergence_empirical} is a plug-in estimator of the Jeffreys divergence in \eqref{J_divergence}, obtained by substituting the unknown class-conditional probabilities with their empirical frequencies. This makes it a standard and consistent approach in information-theoretic inference, which is useful for the hypothesis testing problem \eqref{HT}. Therefore, the estimator \eqref{divergence_empirical} is a test statistic for hypothesis test \eqref{HT}. In the sequel, we present some asymptotic properties of the estimator $J(\hat{\mathbf{p}},\hat{\mathbf{q}})$.
\\

For  simplicity, we introduce the following notations: 
$$\Delta_{p_n}^{a_j}={\hat{p}}_n^{\,a_j} -p_j,\quad \Delta_{q_m}^{a_j}={\hat{q}}_m^{\,a_j} -q_j\quad\textrm{for all}\quad j\in D,$$
$$A_n=\sup\limits_{j \in D}|\Delta_{p_n}^{a_j}|,\quad B_m=\sup\limits_{j \in D}|\Delta_{q_m}^{a_j}| \quad\textrm{and}\quad C_{n,m}=\max\{A_n,\,B_m\}.$$

\noindent
{\bf Theorem 1.} (Almost sure consistency).\,\,\,{\it Let $\mathbf{p}$ and $\mathbf{q}$ be two probability distributions as defined in \eqref{p_q} and satisfying assumption \eqref{assumption}. Furthermore, let ${\hat{\mathbf{p}}}_n$ and ${\hat{\mathbf{q}}}_m$ be generated by i.i.d. sample $(X_1,Y_1), \cdots, (X_N,Y_N)$ according to $\mathbb{P}$ and given by \eqref{empirical_p} and \eqref{empirical_q}. Then the asymptotic result  
\begin{equation}
    \limsup\limits_{(n,m)\rightarrow (+\infty,+\infty)} \frac{| J(\hat{\mathbf{p}}_n,\hat{\mathbf{q}}_m)-J(\mathbf{p},\mathbf{q})|}{C_{n,m}}\leq \sum_{j\in D}\bigg\{\Big|1+\ln\frac{p_j}{q_j}\Big|+\Big|1+\ln\frac{q_j}{p_j}\Big|+\frac{p_j^2+q_j^2}{p_j\,q_j}\bigg\}
\end{equation}
holds almost surely.}
\\

\noindent
{\bf Theorem 2.} (Asymptotic normality).\,\,\,{\it Under the same assumptions as in Theorem 1, the following central limit theorem hold: for $(n,m)\longrightarrow (+\infty,+\infty)$ and $\frac{n\,m}{n+m}\longrightarrow \gamma\in(0,1)$,
\begin{equation}
    \bigg(\frac{n\,m}{m\textrm{V}_1(\mathbf{p}, \mathbf{q})+n\textrm{V}_2(\mathbf{p}, \mathbf{q})}\bigg)^{\frac{1}{2}}\Big(J(\hat{\mathbf{p}}_n,\hat{\mathbf{q}}_m)-J(\mathbf{p},\mathbf{q})\Big) \overset{d}\Longrightarrow\mathcal{N}(0,1),
\end{equation}
where
\begin{eqnarray}
\mathrm{V}_{1}(\mathbf{p},\mathbf{q})&=& \sum_{j \in D}p_j(1-p_j)\bigg\{\Big(1+\ln\frac{p_j}{q_j}\Big)^2+\Big(\frac{q_j}{p_j}\Big)^2\bigg\}\\ \nonumber 
&&\qquad \quad\qquad\qquad \quad-2\sum_{\{(i,j)\in D^2: \,i< j\}}\bigg\{p_ip_j\Big(1+\ln \frac{p_i}{q_i}\Big)\Big(1+\ln \frac{p_j}{q_j}\Big)+q_iq_j\bigg\},
\end{eqnarray}
and and also,
\begin{eqnarray}
\mathrm{V}_{2}(\mathbf{p},\mathbf{q})&=& \sum_{j \in D}q_j(1-q_j)\bigg\{\Big(1+\ln\frac{q_j}{p_j}\Big)^2+\Big(\frac{p_j}{q_j}\Big)^2\bigg\}\\ \nonumber 
&&\qquad \quad\qquad\qquad \quad-2\sum_{\{(i,j)\in D^2: \,i< j\}}\bigg\{q_iq_j\Big(1+\ln \frac{q_i}{p_i}\Big)\Big(1+\ln \frac{q_j}{p_j}\Big)+p_ip_j\bigg\}.
\end{eqnarray}}
\\
The proofs of all our results are postponed to the Appendix.
\\

Using the asymptotic distribution presented in Theorem 2, a  statistical test can be easily derived for the hypothesis test \eqref{HT}, $H_o: \mathbf{p}=\mathbf{q}$. It turns out that when $n$ and $m$ are large, if replacing $\mathbf{p}$ and $\mathbf{q}$ by  by their estimators ${\hat{\mathbf{p}}}_n$ and ${\hat{\mathbf{q}}}_m$ in the computation of $V_1$ and $V_2$, denoted by $\textrm{V}_1({\hat{\mathbf{p}}}_n,{\hat{\mathbf{q}}}_m)$ and $\textrm{V}_2({\hat{\mathbf{p}}}_n,{\hat{\mathbf{q}}}_m)$, has little effect on the normal asymptotic distribution. This leads to the following useful result.
\\

\noindent
{\bf Corollary 1.} (J-Divergence test). {\it Under the null hypothesis $H_0$ in \eqref{HT}, for the test statistic \eqref{divergence_empirical} we have that
\begin{equation}
  \sqrt{\frac{ n m }{ (n+m)\hat{\Sigma}_{n,m} }} J(\hat{\mathbf{p}}_n,\hat{\mathbf{q}}_m) \overset{d}{\Longrightarrow} \mathcal{N}(0,1)
 \quad \textrm{as} \quad n,m \longrightarrow \infty,
\end{equation}
where $$\hat{\Sigma}_{n,m}=\frac{m\textrm{V}_1({\hat{\mathbf{p}}}_n,{\hat{\mathbf{q}}}_m)+n\textrm{V}_2({\hat{\mathbf{p}}}_n,{\hat{\mathbf{q}}}_m)}{n+m}.$$}
\noindent
\\
Therefore, the critical region of the test \eqref{HT} at the significance level $\alpha$ is
\begin{equation}\label{RC}
    \mathcal{R}=\{\hat{J}: J(\hat{\mathbf{p}}_n,\hat{\mathbf{q}}_m)>Q(1-\alpha)\},
\end{equation}
    where $Q(1-\alpha)$ is the $(1-\alpha)$--quantile of the distribution of $J(\hat{\mathbf{p}}_n,\hat{\mathbf{q}}_m)$ under the null hypothesis $H_0$, i.e., $\sqrt{\frac{ n m }{ (n+m)\hat{\Sigma}_{n,m} }} J(\hat{\mathbf{p}}_n,\hat{\mathbf{q}}_m)\sim \mathcal{N}\big(0,1\big)$ asymptotically. Consequently, if we have $J(\hat{\mathbf{p}}_n,\hat{\mathbf{q}}_m) \in \mathcal{R}$, then we reject $H_0$. Therefore, the hypothesis test can be conducted by comparing the test statistic to critical values from the standard normal distribution. Therefore, this test is a two-sided test, where we reject the null hypothesis \( H_0: p = q \) if the test statistic falls into either tail of the standard normal distribution, indicating a significant difference in either direction.

\textbf{Computational Efficiency:}  
In high-dimensional data, where the number of variables is typically very large, the computation of the proposed test for each variable could become computationally intensive. However, the test is independent for each variable and therefore lends itself well to parallelization. Each variable’s test statistic can be computed separately, making it highly scalable for large datasets when computational resources are available. Additionally, in typical applications, the full set of predictors is reduced in a pre-modeling phase using methods such as Information Value (IV), which ranks variables based on their predictive power. This pre-filtering step significantly reduces the number of variables considered in subsequent tests, thereby minimizing the computational burden.

\section{Performance of the test on simulated data}\label{simulation}

In this section we present the results of a simulation study in order to evaluate the performance of the test \eqref{HT}. For greater control of the simulations, we assume that 
\begin{equation}\label{multinomial_model}
    \Bigg(\sum_{i=1}^{\Lambda}X_i=a_1,\dots ,\sum_{i=1}^{\Lambda}X_i=a_r\Bigg)\Bigg|Y=y
 \sim \textrm{Multinomial}\,(\Lambda, \Theta) , 
\end{equation}
where $\Lambda=n\cdot\mathbf{I}[y=1]+m\cdot\mathbf{I}[y=0]$. Moreover, where $\Theta=\Big(z_{yj}\,\theta_y^j(1-\theta_y)^{r-j}\Big)_{1\leq j \leq r}$, $\theta_y \in (0,1)$ and $z_{yj}$ is a normalization constant for $(y,j)\in \{0,1\}\times D$. By assuming \eqref{multinomial_model}, we have $\mathbf{p}$ and $\mathbf{q}$ are determined in such a way that 
\begin{equation}\label{probs_thetas}
    p_j=z_{1j}\,\theta_1^j(1-\theta_1)^{r-j}\qquad\textrm{and} \qquad q_j=z_{0j}\,\theta_0^j(1-\theta_0)^{r-j},
\end{equation}
where $j=1,\dots, r$. The empirical distributions $\hat{\mathbf{p}}_n$ and $\hat{\mathbf{q}}_m$ be generated according to \eqref{multinomial_model} and given by \eqref{empirical_p} and \eqref{empirical_q} respectively. To run the simulations, the value of $\theta_1=0.5$ is set, from which $\hat{\mathbf{p}}_n$ is generated, and a value for $\theta_0=\theta$ is selected from a grid on $[0,1]$ with step length $10^{-3}$, from which $\hat{\mathbf{q}}_m$ is generated. For each point $(\theta_1,\theta)\in \{0.5\}\times [0,1]$, we obtain $M$ simulations for which we measure the divergences $J(\hat{\mathbf{p}}_n,\hat{\mathbf{q}}_m)$. From \eqref{RC}, for each $\theta \in [0,1]$, we define the empirical power function as
\begin{equation}
    \textrm{Power-Function}=\frac{1}{M}\,\sum\limits_{i=1}^M \mathbf{I}[J(\hat{\mathbf{p}}_n,\hat{\mathbf{q}}_m)\in \mathcal{R}].
\end{equation}
where $M=10^3$ for all simulations. Furthermore, to control and verify the effect of the imbalance in our simulations we define
\begin{equation}
    \textrm{Imbalance-Rate}=\frac{n}{n+m}.
\end{equation}

The goal of the simulations is to assess the behavior of the test under various conditions, including different imbalance levels and numbers of bins. Before presenting the results, it is important to clarify the conditions under which the simulations were conducted under the null hypothesis \( H_0: \mathbf{p}=\mathbf{q} \). In these simulations, the class-conditional probabilities were set equal, i.e., \( \mathbf{p}=\mathbf{q} \), and the parameter \( \theta \) was varied across a range of values to simulate different scenarios where no effect exists between the distributions of \( p \) and \( q \). Specifically, \( \theta = 0.5 \) represents the case of no difference between the two distributions.

\begin{remark}
 Note that assuming \eqref{multinomial_model} and \eqref{probs_thetas}
 does not restrict our results since this model is used only to generate the empirical distributions and not for test evaluation itself.The reason why we have chosen a parametric model for our experiments is for greater control of the simulations and simplicity in the resulting graphs. As mentioned above, our test does not make specific assumptions about the shape of the distributions evaluated.   
\end{remark}

\begin{figure}[h]
\centering
\includegraphics[scale=0.55]{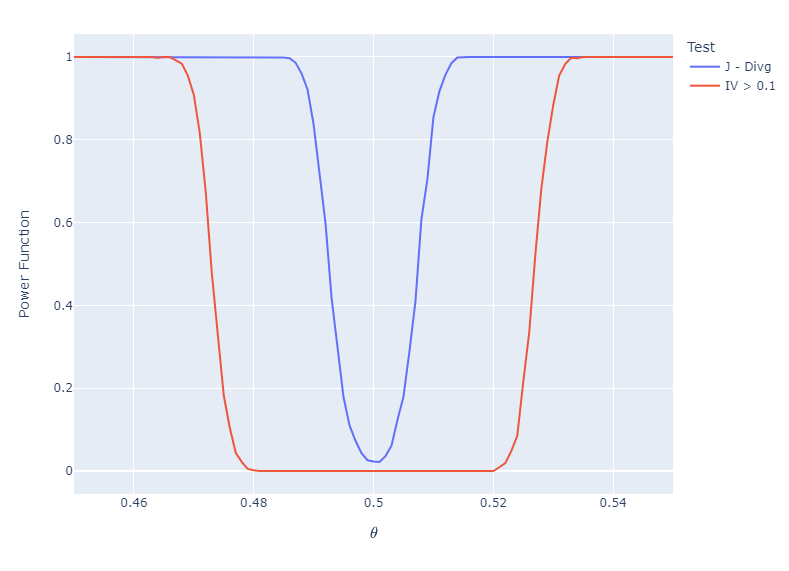}
\caption{Comparison between the power function of the J-Divergence test (top line) and the empirical criterion IV $> 0.1$ (bottom line). The solid line at the top of the plot corresponds to the J-Divergence test, while the solid line at the bottom represents the IV $> 0.1$ criterion. This description ensures clarity when printed in black and white, as the lines can be distinguished based on their relative positions. Simulation parameters: \( N = 50\,300 \), \( n = 300 \), \( m = 50\,000 \), Imbalance-Rate = 0.00596, \( \alpha = 0.1\% \), \( r = 10 \).}
\label{J-IV}
\end{figure}

Before presenting the simulation results, we clarify the selection criterion used for IV. In practice, IV thresholds are commonly applied to eliminate irrelevant predictors and retain those with meaningful contribution. In our experiments, we adopt the filtering rule \(0.1 < \text{IV} < 0.5\). The lower threshold eliminates weak predictors, while the upper limit helps avoid overly dominant variables that may cause overfitting. This range is chosen to support the construction of stable and generalizable models. For simplicity, although we continue to refer to the threshold as “IV $> 0.1$” throughout the rest of the paper, it should be understood that all variables considered satisfy the condition \(0.1 < \text{IV} < 0.5\).

\begin{figure}[h]
\centering
\includegraphics[scale=0.55]{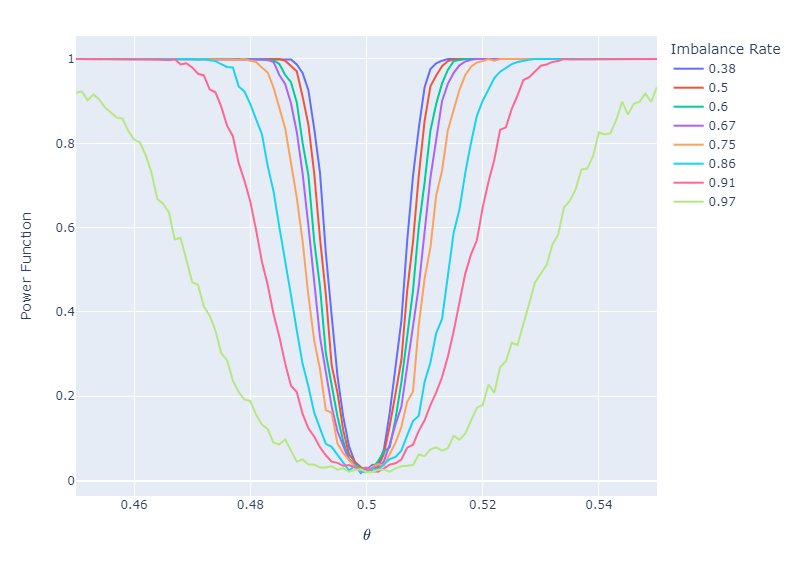}
\caption{Power function of the J-Divergence test for different imbalance rates. The lines represent different imbalance rates, with the line at the top of the plot corresponding to an imbalance rate of 0.38, the next line representing 0.5, followed by 0.6, 0.67, 0.75, 0.86, and the bottom line representing 0.97. This description ensures clarity in black and white prints, as the lines can be distinguished based on their relative positions. Simulation parameters: \( r = 10 \), \( \alpha = 0.1\% \), \( n = 3000 \), \( m \in [100, 50\,00] \).}

\label{J_inbalance}
\end{figure}

\begin{figure}[h]
\centering
\includegraphics[scale=0.44]{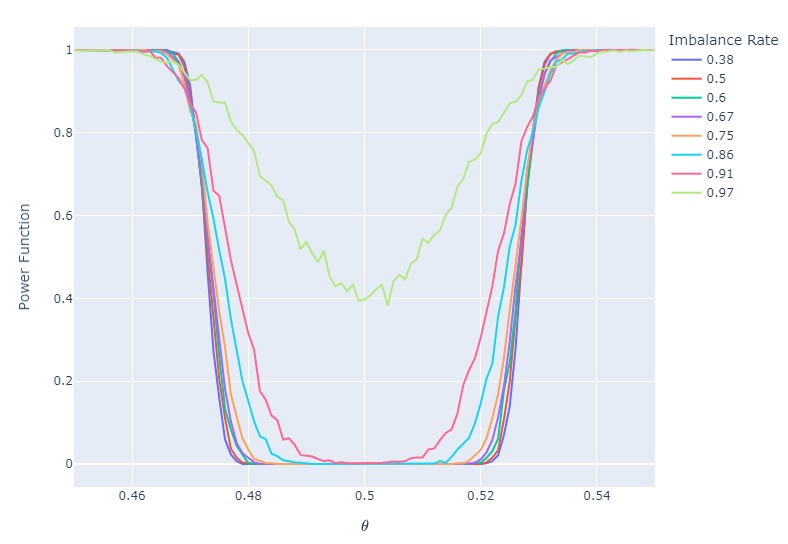}
\caption{Power function of criterion IV $> 0.1$ for different imbalance rates. The lines represent different imbalance rates, with the topmost line corresponding to an imbalance rate of 0.38, followed by 0.5, 0.6, 0.67, 0.75, 0.86, and the bottommost line representing 0.97. This line style differentiation ensures clarity when printed in black and white, as the lines can be distinguished based on their relative positions. Simulation parameters: \( r = 10 \), \( \alpha = 0.1\% \), \( n = 3\,000 \), \( m \in [100, 5\,000] \).}

\label{IV_inbalance}
\end{figure}

In the first simulation example, we compare the power function of our test with the power function of the criterion $\textrm{IV}>0.1$ detailed in Section \ref{intro}. In this case, we measure the power function as
\begin{equation}
    \textrm{Power-Function}=\frac{1}{M}\,\sum\limits_{i=1}^M \mathbf{I}[J(\hat{\mathbf{p}}_n,\hat{\mathbf{q}}_m)> 0.1],
\end{equation}
for each $\theta\in [0,1]$. The results of this comparison are shown in Figure \ref{J-IV}. From this figure we have that, in conditions of imbalance, the superiority of the performance of our test compared to the commonly used criterion is notable. It is natural to wonder if this scenario holds at different levels of imbalance. Therefore, in Figure \ref{J_inbalance} and Figure \ref{IV_inbalance} we evaluate the effect of different levels of imbalance on the performance of both methodologies. On the one hand, the performance of our test is robust to different intensities of the imbalance, see Figure \ref{J_inbalance}. However, the performance of criterion IV$>0.1$ is drastically affected by high levels of imbalance, see Figure \ref{IV_inbalance}. Note that the green line in Figure \ref{IV_inbalance} starts almost at 0.4, which represents the high Type-I error, high levels of false positives, associated with this criterion. It is worth mentioning that our test does not experience this Type-I error inflation generated by high levels of imbalance. In this sense, we say that our test is much more robust than the conventionally used criterion IV$>0.1$.
\begin{remark}
    The hypothesis test that we propose in this work is much more robust to the presence of high levels of imbalance in the target variable.
\end{remark}

\begin{figure}[h]
\centering
\includegraphics[scale=0.52]{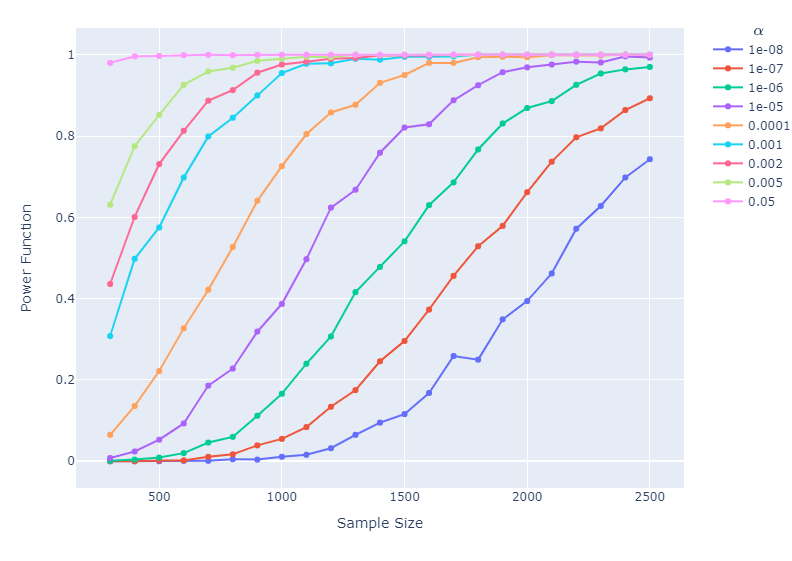}
\caption{Power function of the J-Divergence test for different values of \( \alpha \). The lines represent different values of \( \alpha \), with the line at the top of the plot corresponding to \( \alpha = 10^{-8} \), followed by \( \alpha = 10^{-7} \), \( \alpha = 10^{-6} \), \( \alpha = 10^{-5} \), \( \alpha = 10^{-4} \), \( \alpha = 0.001 \), \( \alpha = 0.002 \), \( \alpha = 0.005 \), and the line at the bottom representing \( \alpha = 0.05 \). This line style differentiation ensures clarity when printed in black and white, as the lines can be distinguished based on their relative positions. Simulation parameters: \( n = m \in [300, 2\,500] \), Imbalance-Rate = 0.5, \( r = 10 \).}
\label{PF_alphas}
\end{figure}

\begin{figure}[h]
\centering
\includegraphics[scale=0.52]{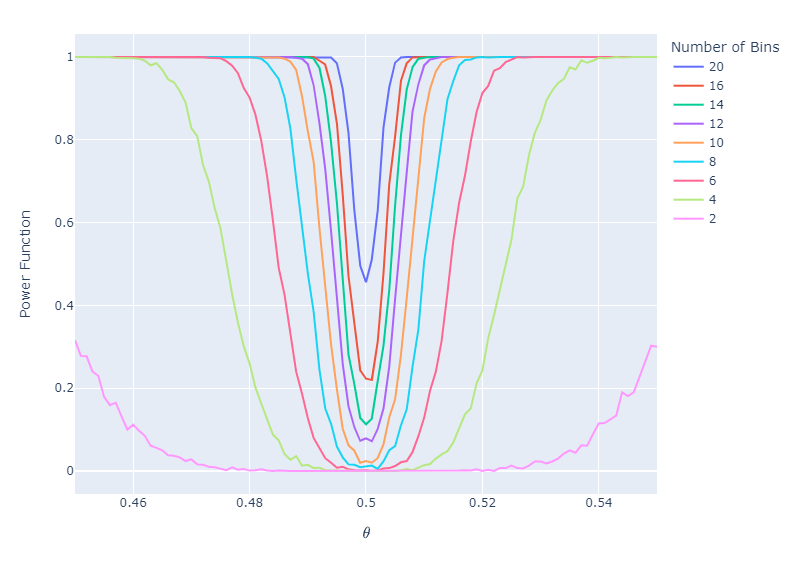}
\caption{Power function of the J-Divergence test for different numbers of bins. The lines represent different numbers of bins, with the line at the top corresponding to 20 bins, followed by 16, 14, 12, 10, 8, 6, 4, and the bottommost line representing 2 bins. This line style differentiation ensures clarity when printed in black and white, as the lines can be distinguished based on their relative positions. Simulation parameters: \( n = m = 3\,000 \), Imbalance-Rate = 0.5, \( r \in [2, 20] \).}
\label{PF_bins}
\end{figure}

It is also of primary interest to know how the performance of our test depends on the levels of statistical significance $\alpha$ and the number of bins $r$. Our test performs best for $\alpha \in [0.001\%, 0.1\%]$, see Figure \ref{PF_alphas}. Furthermore, our test shows robust performance, without significant Type-I error, for bin numbers $r\in [2, 14]$, see Figure \ref{PF_bins}.
\begin{remark}
    Although the optimal value for the level of statistical significance $\alpha$ and the number of bins $r$ that should be used to guarantee a good performance of the hypothesis testing proposed in this work should vary case by case, we suggest using $\alpha=0.01\%$ and $r\leq 14$.
\end{remark}

It is worth emphasizing that the simulation study was designed with a specific purpose: to examine the behavior of the proposed J-Divergence test in scenarios where traditional IV-based criteria often fail—namely, under conditions of class imbalance. While IV thresholds may perform reasonably well in balanced settings, our goal was to highlight the robustness of our method in more challenging contexts. Furthermore, the theoretical properties of our test (almost sure consistency and asymptotic normality) provide solid guarantees under minimal assumptions. Therefore, the simulations presented in this section are meant to illustrate the practical impact of these theoretical results, particularly in cases where empirical heuristics tend to be less reliable.

\section{Features selection in  fraud detection}\label{aplication}

The data analyzed in this section correspond to a public dataset provided by Vesta Corporation as part of a fraud detection challenge hosted on Kaggle. It consists of real-world e-commerce transactions, including a wide variety of features such as product codes, payment methods, device types, address matching indicators, and anonymized variables derived from card usage and transactional behavior. These variables capture both behavioral patterns and device-related characteristics that may help identify fraudulent activity. The dataset is divided into two subsets: a training set comprising 472,432 operation records and a test set with 118,108 records. Both subsets share the same structure, consisting of 369 features and a binary target variable indicating whether a transaction was fraudulent. The target variable is highly imbalanced, with a fraud rate of approximately $0.3\%$, which is standard in real-world fraud detection problems. This imbalance poses additional challenges for building robust machine learning classifiers capable of accurately detecting fraud, see \citet{Addison2019}. 

\begin{figure}[h]
\centering
\includegraphics[scale=0.40]{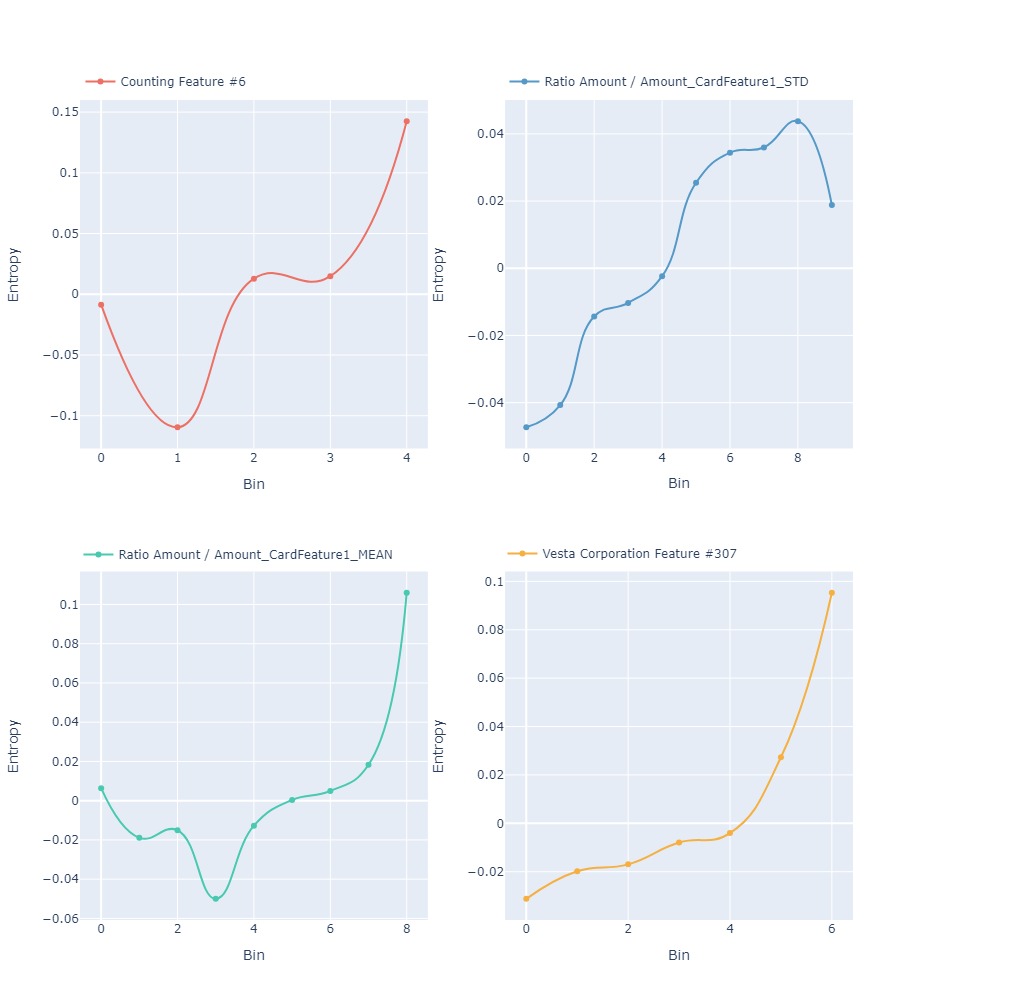}
\caption{Shannon entropy versus bins for four features selected only by the J-Divergence test. For each graph, the horizontal axis contains the mean within each bin and the vertical axis contains the Shannon entropy measure within each bin with respect to the target variable.}
\label{Entropy_Bins}
\end{figure}

The challenge is to build classification machine learning models to detect fraud with high precision and simultaneously with low levels of false positives. To this end, the first step is to features selection based on predictive power, that is, prioritize the features that contain the most information and  that will allow correct discrimination of fraudulent operations. For this first stage we have considered three alternatives: i) Do not use any characteristic selection criteria, use all available features. ii) Use the criterion IV$>0.1$. iii) Use the J-Divergence test that we have proposed in this work. The results for each selection criterion evaluated are presented in Table \ref{table_selection}. Our test selected 42 features that the IV$>0.1$ criterion did not select. A very interesting fact is that these 42 features not considered by the second criterion contain relevant information for the detection of fraud. For example, to have a general idea, of these 42 additional features we have selected the four with the lowest p-value assigned by our test, see Table \ref{ouput_statistical_iv}. These four features have a clear non-linear trend with respect to the entropy within each bin, this confirms the significant predictive power of these characteristics when identifying fraud, see Figure \ref{Entropy_Bins}. These same patterns are repeated in the remaining variables selected by our test. Therefore, in this first stage, prior to modeling, our test satisfactorily selected features that contain relevant information for the classification problem and showed superiority with respect to second methodology.
\begin{table}[]
\centering 
\begin{tabular}{|c|c|c|}
\hline
\textbf{Without Selection} & \textbf{Criterion IV $>$ 0.1} & \textbf{J-Divergence Test} \\ \hline
369                        & 220                         & 262                        \\ \hline
\end{tabular}
\caption{Number of features selected for each of the three methodologies in the fraud detection application of Chapter 4: no selection, IV $> 0.1$, and J-Divergence test.}
\label{table_selection}
\end{table}

\begin{table}[]
\begin{tabular}{|c|c|c|c|}
\hline
Predictor                                 & J - Estimate & Std - Error & p - value  \\ \hline
Counting Feature \#6                      & 0.099010     & 0.003176   & 2.694e$^{-213}$ \\ \hline
Ratio Amount/Amount\_CardFeature1\_STD  & 0.085637     & 0.002939   & 1.036e$^{-186}$ \\ \hline
Ratio Amount/Amount\_CardFeature1\_MEAN & 0.077940     & 0.002806   & 8.449e$^{-170}$ \\ \hline
Vesta Corporation Feature 307\#           & 0.080216     & 0.002893   & 3.681e$^{-169}$ \\ \hline
\end{tabular}
\caption{Results of the J-Divergence test for the selected variables in the fraud detection application. The table shows the J-estimates, standard errors, and p-values for each predictor selected by the J-Divergence test.}
\label{ouput_statistical_iv}
\end{table}

In the second stage, the modeling stage, for each selected feature subset, we will use the train dataset to train a Light Gradient Boosting Machine (LightGBM) with their respective parameters optimized for each case
with their respective parameters optimized for each case. The main LightGBM performance indicators on test dataset are shown in Table \ref{performance_LGBM}. From this table we highlight that our test has a little more than 1\% in precision with respect to the second methodology IV$>0.1$, and judging by the other indicators, this increase is due to the decrease in false positives, which in this context translate into a non-negligible decrease in false fraud alerts. As mentioned above, a reduction of more than 1\% in false fraud alerts generally has a considerable economic and financial impact since it avoids the blocking of many genuine paymets operations and all the inconvenience that this generates for customers. If we compare our test with the first methodology, that is, training the LightGBM without any prior feature selection criteria, we have a 3\% advantage in precision with 107 fewer features. In addition to the economic impact generated by the increase in precision in fraud detection, we must consider the reduction in engineering cost, implementation cost, production and maintenance, as a consequence of the reduction of more than 100 features in the classifier. Therefore, we can argue that the appropriate use of the test proposed in this work leads to the features selection efficient prior to the modeling stage.

\begin{table}[]
\centering
\begin{tabular}{|c|c|c|c|}
\hline
\textbf{Metric}  & \textbf{Without Selection} & \textbf{Criterion IV $>$ 0.1} & \textbf{J-Divergence Test} \\ \hline
Precision        & 0.85496                     & 0.87141                     & 0.88491                     \\ \hline
Recall          & 0.74594                     & 0.73457                     & 0.73812                     \\ \hline
AUC             & 0.97098                     & 0.97112                     & 0.97188                     \\ \hline
F1 Score        & 0.79624                     & 0.79666                     & 0.79734                     \\ \hline
\end{tabular}
\caption{LightGBM performance on the test dataset for different feature selection methods. The table compares the performance of LightGBM using no feature selection, the IV $> 0.1$ criterion, and the J-Divergence test. The metrics evaluated include Precision, Recall, AUC, and F1 score.}
\label{performance_LGBM}
\end{table}

In this study, we applied the J-Divergence test for feature selection and proceeded to model the data using LightGBM. However, it is important to emphasize that the feature selection methodology employed here is independent of the machine learning model selected. Our approach, which utilizes IV and the J-Divergence test for variable filtering, is designed to be applicable to any model chosen during the modeling stage. This ensures flexibility and robustness in the pre-modeling phase. Given the nature of the dataset, which is highly imbalanced, we selected LightGBM due to its efficiency, scalability, and ability to handle large and complex datasets. Moreover, LightGBM is one of the most widely used classifiers in credit scoring and fraud detection, making it a natural choice for this case study. Nonetheless, it should be noted that our methodology is not limited to LightGBM and could be applied to any model selected for the final predictions, highlighting the flexibility and adaptability of our approach.

\begin{table}[]
\centering
\begin{tabular}{|l|c|c|c|c|c|}
\hline
\textbf{Filter} & \textbf{\# Features} & \textbf{Precision} & \textbf{Recall} & \textbf{AUC} & \textbf{F1} \\
\hline
IV $>$ 0.02 & 312 & 0.86950 & 0.73610 & 0.97085 & 0.79512 \\
\hline
\end{tabular}
\caption{Summary performance for exploratory 0.02 threshold.}
\label{tab:iv002_summary}
\end{table}

To assess whether lower IV thresholds could serve as an appropriate baseline, we include an additional evaluation based on the threshold IV $> 0.02$, which is frequently cited in the literature as a minimum cutoff to exclude variables with weak univariate relevance, e.g., see \citet{pang2020borrowers, li2024financial}. This more permissive criterion resulted in the selection of 312 features. However, as shown in Table~\ref{tab:iv002_summary}, this strategy did not yield better performance compared to the benchmark IV $> 0.1$. The results suggest that IV $> 0.02$ tends to retain variables with weak discriminative power, which may introduce noise rather than improving model accuracy. Therefore, IV $> 0.1$ remains a more stable and effective filtering rule in our setting, and a more appropriate baseline for comparison with the proposed J-Divergence test.

In the context of credit risk modeling, \citet{hade2025} propose an integrative feature selection framework that combines multiple criteria—namely, Kruskal-Wallis, Cramér’s V, and Information Value—into a composite score used to rank variables. Although this approach captures different aspects of relevance, it requires normalization of heterogeneous measures and the calibration of empirical thresholds. In contrast, the J-Divergence test introduced in this paper offers a unified, nonparametric testing framework with clear asymptotic guarantees and interpretability via \(p\)-values. This makes it suitable not only as a standalone filter but also as a complementary component within integrative pipelines like the one proposed by \citet{hade2025}. Our methodology could thus enhance such procedures by providing both theoretical justification and operational robustness, particularly in imbalanced settings such as fraud detection.

\textbf{Limitations:} While the J-Divergence test offers a theoretically grounded and practically effective method for variable selection—particularly in imbalanced scenarios—some limitations must be acknowledged.  First, the method requires input features to be categorical or discretized into a finite number of bins. For continuous variables, this entails a preprocessing step that involves partitioning the domain into intervals. The reliability of the test depends critically on the quality of this discretization. Second, the test relies on the availability of sufficient data. Since the underlying theoretical results, consistency and asymptotic normality, are based on empirical distributions over bins, each bin must contain a non-negligible number of observations. In datasets with low sample size or very sparse features, it may be difficult to construct bins with enough support to guarantee valid inference. In such cases, the power of the test may degrade or the asymptotic approximations may be unreliable. However, this limitation is not specific to our method, but rather inherent to any statistical test based on empirical distributions or nonparametric estimates. In fact, model-based feature selection techniques, e.g., wrapper methods or tree-based importance scores, also suffer from instability or overfitting when the data is scarce or highly unbalanced. Our approach offers a clean and interpretable alternative precisely in the context where traditional empirical IV filtering lacks formal guarantees. Finally, the test is currently formulated for binary classification tasks. Extending the method to multiclass settings or continuous response variables remains an open direction for future work.

\section{Implementation}\label{python}
In order to facilitate the usability and reproducibility of the main results of this work, we have developed  the open-source python library ``statistical-iv" available in \url{https://pypi.org/project/statistical-iv/}.

\section{Conclusions and future work}
In this paper we propose a non-parametric hypothesis test for the IV. To our knowledge, there is no work in the literature in which a formal statistical test for this quantity has been presented. We derive efficient asymptotic formulas for the test statistic that imply that our test is consistent. Furthermore, in different scenarios, simulations show that our test easily surpasses existing empirical criteria about IV based on fixed thresholds. As imbalance is a recurring scenario in classification problems, our test works very well even for scenarios with high levels of imbalance. In a real fraud identification data set, we show the potential of our J-Divergence statistic to evaluate the predictive power of the considered features.

Although the main focus of this article is binary classification, the generalization of the J test statistic, and therefore of the IV, to the case of multinomial classification is possible. Perhaps the idea of this generalization is simple and natural, although we believe that the technical details are not so trivial and warrants future research.

In summary, our proposal is designed as a robust and interpretable tool for the pre-modeling phase of feature selection. The J-Divergence test provides a formal statistical alternative to empirical IV thresholds, without relying on model-specific assumptions. This makes it particularly suitable for early-stage variable prioritization, especially in imbalanced datasets where conventional practices may fail. Rather than replacing model-based selection methods such as Lasso or tree-based importance, our methodology complements them by offering a theoretically grounded filtering step that can improve the quality and reliability of downstream modeling efforts. As documented in the literature \citet{guyon2003introduction, chandrashekar2014survey}, filter-based methods like ours are widely used to guide the selection of variables before engaging in more complex model-based selection processes.

\appendix
\section*{Appendix: Proofs}

\subsection*{A1. Proof of Theorem 1}
 Let $\phi(x,y)=(x-y) \ln\big(\frac{x}{y}\big)$, in such a way that
\begin{equation}
    J(\hat{\mathbf{p}}_n,\hat{\mathbf{q}}_m)=\sum_{j\in D}\phi({\hat{p}}_n^{\,a_j},{\hat{q}}_n^{\,a_j})\quad\textrm{and}\quad  J({\mathbf{p}},{\mathbf{q}})=\sum_{j\in D}\phi(p_j,q_j).
\end{equation}
For a fixed $j \in D$ and making use of the continuous differentiability of the function $\phi(x,y)$ along with first-order Taylor`s theorem, see \citet{marsden2003vector}, we get
\begin{eqnarray}\label{Taylor}
\phi({\hat{p}}_n^{\,a_j},{\hat{p}}_n^{\,a_j})&=&\phi(p_j+\Delta_{p_n}^{a_j},\,q_j+\Delta_{q_n}^{a_j})\\ \nonumber 
&=&\phi(p_j,q_j)+\Delta_{p_n}^{a_j}\,\partial_{x}\phi(p_j,q_j)+\Delta_{q_n}^{a_j}\,\partial_{y}\phi(p_j,q_j)+R_j(n,m),
\end{eqnarray}
where, for $0<\theta_j, \eta_j<1$, the Lagrange’s remainder term in the following form
\begin{eqnarray}\label{R}
R_j(n,m)&=&\frac{1}{2}\Big[\big(\Delta_{p_n}^{a_j}\big)^2\partial_{xx}\phi\big(p_j+\theta_j\Delta_{p_n}^{a_j},\,q_j+\eta_j\Delta_{q_n}^{a_j}\big)\\ \nonumber 
&&\qquad\qquad\qquad+2\,\Delta_{p_n}^{a_j}\Delta_{q_n}^{a_j}\,\partial_{xy}\phi\big(p_j+\theta_j\Delta_{p_n}^{a_j},\,q_j+\eta_j\Delta_{q_n}^{a_j}\big)\\ \nonumber
&&\qquad\qquad\qquad\qquad\qquad\qquad+ \big(\Delta_{q_n}^{a_j}\big)^2\partial_{yy}\phi\big(p_j+\theta_j\Delta_{p_n}^{a_j},\,q_j+\eta_j\Delta_{q_n}^{a_j}\big)\Big].
\end{eqnarray}
Summing over $j\in D$ in \eqref{Taylor} we have that
\begin{equation}\label{J_equality}
    J(\hat{\mathbf{p}}_n,\hat{\mathbf{q}}_m)-J(\mathbf{p},\mathbf{q})=\sum_{j\in D}\Delta_{p_n}^{a_j}\,\partial_{x}\phi(p_j,q_j)+\sum_{j\in D}\Delta_{q_n}^{a_j}\,\partial_{y}\phi(p_j,q_j)+\sum_{j\in D}R_j(n,m),
\end{equation}
which leads to
\begin{equation}\label{J_inequality}
    |J(\hat{\mathbf{p}}_n,\hat{\mathbf{q}}_m)-J(\mathbf{p},\mathbf{q})|\leq A_n\sum_{j\in D}|\partial_{x}\phi(p_j,q_j)|+B_n\sum_{j\in D}|\partial_{y}\phi(p_j,q_j)|+\sum_{j\in D}|R_j(n,m)|.
\end{equation}
In \eqref{J_inequality}, notice that $\sum_{j\in D}|\partial_{x}\phi(p_j,q_j)|<\infty$ and $\sum_{j\in D}|\partial_{y}\phi(p_j,q_j)|<\infty$. In the sequel, it is easy to deduce that the third term on the right hand side of equation \eqref{J_inequality} satisfies the following inequality
\begin{eqnarray}\label{R_inequality}
\frac{1}{C_{n,m}}\sum_{j\in D}|R_j(n,m)|&\leq&\frac{1}{2}\bigg[A_n\,\sum_{j\in D}\big|\partial_{xx}\phi\big(p_j+\theta_j\Delta_{p_n}^{a_j},\,q_j+\eta_j\Delta_{q_n}^{a_j}\big)\big|\\ \nonumber 
&&\qquad+2\,A_n\,\sum_{j\in D}\big|\partial_{xy}\phi\big(p_j+\theta_j\Delta_{p_n}^{a_j},\,q_j+\eta_j\Delta_{q_n}^{a_j}\big)\big|\\ \nonumber
&&\qquad\qquad\qquad\quad+ B_m\,\sum_{j\in D}\big|\partial_{yy}\phi\big(p_j+\theta_j\Delta_{p_n}^{a_j},\,q_j+\eta_j\Delta_{q_n}^{a_j}\big)\big|\bigg].
\end{eqnarray}
Based on the strong law of large numbers, we know that both $\Delta_{p_n}^{a_j}$ and $\Delta_{q_m}^{a_j}$ tends to $0$ as $n, \, m\to \infty$ almost surely, and consequently, $A_n$ and $B_m$ also tends to $0$. Based on the same argument, we have $\sum_{j\in D}\big|\partial_{xx}\phi\big(p_j+\theta_j\Delta_{p_n}^{a_j},\,q_j+\eta_j\Delta_{q_n}^{a_j}\big)\big| \to \sum_{j\in D}|\partial_{xx}\phi(p_j,\,q_j)|<\infty$, $\sum_{j\in D}\big|\partial_{xy}\phi\big(p_j+\theta_j\Delta_{p_n}^{a_j},\,q_j+\eta_j\Delta_{q_n}^{a_j}\big)\big| \to \sum_{j\in D}|\partial_{xy}\phi(p_j,\,q_j)|<\infty$, and likewise, $\sum_{j\in D}\big|\partial_{yy}\phi\big(p_j+\theta_j\Delta_{p_n}^{a_j},\,q_j+\eta_j\Delta_{q_n}^{a_j}\big)\big| \to \sum_{j\in D}|\partial_{yy}\phi(p_j,\,q_j)|<\infty$ as $n, \, m\to \infty$ almost surely. Therefore, from equation \eqref{R_inequality} we have to
\begin{equation}\label{R_convergence}
    \frac{1}{C_{n,m}}\sum_{j\in D}|R_j(n,m)|\longrightarrow 0 \quad\textrm{as}\quad n, m \longrightarrow \infty,
\end{equation}
almost surely. Finally, we combine equations \eqref{J_inequality} and \eqref{R_convergence}, which leads us to
\begin{eqnarray}
&&
\limsup\limits_{(n,m)\rightarrow (+\infty, +\infty)} 
\frac{|J(\hat{\mathbf{p}}_n,\hat{\mathbf{q}}_m)-J(\mathbf{p},\mathbf{q})|}{C_{n,m}}
\leq \sum_{j\in D} |\partial_{x}\phi(p_j,q_j)| 
+ \sum_{j\in D} |\partial_{y}\phi(p_j,q_j)| \quad \\ \nonumber 
&& \leq \sum_{j\in D} \Bigg| \ln\Big(\frac{p_j}{q_j}\Big) 
+ \frac{p_j - q_j}{p_j} \Bigg| \\ \nonumber
&& + \sum_{j\in D} \Bigg| \ln\Big(\frac{q_j}{p_j}\Big) 
+ \frac{q_j - p_j}{q_j} \Bigg| \\ \nonumber
&\leq& \sum_{j\in D} \Bigg\{ 
\Big|1+\ln\Big(\frac{p_j}{q_j}\Big)\Big| 
+ \Big|1+\ln\Big(\frac{q_j}{p_j}\Big)\Big| + \frac{q_j}{p_j} + \frac{p_j}{q_j} 
\Bigg\}.
\end{eqnarray}

almost surely. 
\hfill$\square$
\vskip 0.2in
\subsection*{A2. Proof of Theorem 2}
From definition \eqref{empirical_p} we know that the $r$-tuple $n\hat{\mathbf{p}}_n=\big(n\hat{p}_n^{a_j}\big)_{j\in D}$ follows a Multinomial law  with parameters $n$ and $\mathbf{p}$, denoted by $n\hat{\mathbf{p}}_n \sim \mathcal{M}(n,\mathbf{p})$. Similarly, from definition \eqref{empirical_q}, $m\hat{\mathbf{q}}_m \sim \mathcal{M}(m,\mathbf{q})$. Therefore, applying the weak convergence of the Multinomial law, see for example Proposition 2 in \citet{lo2021weak}, we have
\begin{equation}\label{normality_1}
    \Bigg(\sqrt{\frac{n}{p_j}}\Delta_{p_n}^{a_j}\Bigg)_{j\in D}\overset{d}\Longrightarrow \quad Z_\mathbf{p}\sim\mathcal{N}_r(0, \Sigma_\mathbf{p}) \quad \textrm{as}\quad n\longrightarrow\infty,
\end{equation}
and
\begin{equation}\label{normality_2}
    \Bigg(\sqrt{\frac{m}{q_j}}\Delta_{q_m}^{a_j}\Bigg)_{j\in D}\overset{d}\Longrightarrow \quad Z_\mathbf{q}\sim\mathcal{N}_r(0, \Sigma_\mathbf{q}) \quad \textrm{as}\quad m\longrightarrow\infty.
\end{equation}
Here $Z_{\mathbf{p}}=(Z_{p_j})_{j \in D}$ $\Sigma_\mathbf{p}$ and $Z_{\mathbf{q}}=(Z_{q_j})_{j \in D}$ are two independent
Gaussian random vectors whit mean zero and covariance matrix $\Sigma_\mathbf{p}$ and $\Sigma_\mathbf{q}$ respectively, whose elements are
\begin{equation}\label{Sigma_p}
    \big(\Sigma_\mathbf{p}\big)_{i,j}=(1-p_j)\delta_{ij}-\sqrt{p_ip_j}(1-\delta_{ij}),\quad (i,j)\in D^2,
\end{equation}
and
\begin{equation}\label{Sigma_q}
    \big(\Sigma_\mathbf{q}\big)_{i,j}=(1-q_j)\delta_{ij}-\sqrt{q_iq_j}(1-\delta_{ij}),\quad (i,j)\in D^2,
\end{equation}
where $\delta_{ij}$ is the Kronecker delta. Using the linearity properties of the multivariate normal distribution in \eqref{normality_1}, and with the help of \eqref{Sigma_p}, we get
\begin{equation}\label{convergence_p}
\sum_{j\in D}\Delta_{p_n}^{a_j}\,\partial_{x}\phi(p_j,q_j)\overset{d}\Longrightarrow \frac{1}{\sqrt{n}}\sum_{j\in D}\sqrt{p_j}\,\partial_{x}\phi(p_j,q_j)Z_{p_j}\quad \textrm{as}\quad n\longrightarrow \infty,
\end{equation}
which the right side of equation \eqref{convergence_p} is a random variable that is normally distributed with mean zero and variance equal to
\begin{eqnarray}
\frac{1}{n}\,\mathrm{V}_{1}(\mathbf{p},\mathbf{q}) &=& 
\mathbb{V}\textrm{ar} \Bigg( \frac{1}{\sqrt{n}} \sum_{j\in D} \sqrt{p_j} \,\partial_{x}\phi(p_j,q_j)Z_{p_j} \Bigg) \\ \nonumber 
&=& \frac{1}{n} \bigg\{ 
\sum_{j\in D} \mathbb{V}\textrm{ar} \big( \sqrt{p_j}\,\partial_{x}\phi(p_j,q_j)Z_{p_j} \big) \\ \nonumber
&&\quad + 2 \sum_{\{(i,j)\in D^2: \,i< j\}} 
\mathbb{C}\textrm{ov} \Big( \sqrt{p_i}\,\partial_{x}\phi(p_i,q_i)Z_{p_i}, 
\sqrt{p_j}\,\partial_{x}\phi(p_j,q_j)Z_{p_j} \Big) 
\bigg\} \\ \nonumber 
&=& \frac{1}{n} \bigg\{ 
\sum_{j\in D} p_j(1-p_j)\,\big(\partial_{x}\phi(p_j,q_j)\big)^2 \\ \nonumber
&& \quad - 2 \sum_{\{(i,j)\in D^2: \,i< j\}} 
p_i p_j\,\partial_{x}\phi(p_i,q_i)\,\partial_{x}\phi(p_j,q_j) 
\bigg\} \\ \nonumber
&=&  \frac{1}{n} \Bigg\{ 
\sum_{j\in D} p_j(1-p_j) 
\bigg[ \bigg( 1+\ln\frac{p_j}{q_j} \bigg)^2 
+ \Big(\frac{q_j}{p_j}\Big)^2 \bigg] \\ \nonumber 
&&\quad - 2 \sum_{\{(i,j)\in D^2: \,i< j\}} 
p_i p_j \bigg[ \big(1+\ln\frac{p_i}{q_i} \big) 
\big(1+\ln\frac{p_j}{q_j} \big) 
+ \frac{q_i q_j}{p_i p_j} \bigg] 
\Bigg\}.
\end{eqnarray}
Similarly, from equations \eqref{normality_2} and \eqref{Sigma_q}, we know that
\begin{equation}\label{convergence_q}
\sum_{j\in D}\Delta_{q_m}^{a_j}\,\partial_{y}\phi(p_j,q_j)\overset{d}\Longrightarrow \frac{1}{\sqrt{m}}\sum_{j\in D}\sqrt{q_j}\,\partial_{y}\phi(p_j,q_j)Z_{q_j}\quad \textrm{as}\quad m\longrightarrow \infty,
\end{equation}
which the right side of equation \eqref{convergence_q} is a random variable that is normally distributed with mean zero and variance equal to
\begin{eqnarray}
\frac{1}{m}\,\mathrm{V}_{2}(\mathbf{p},\mathbf{q}) &=& 
\mathbb{V}\textrm{ar}\Bigg( \frac{1}{\sqrt{m}}\sum_{j\in D} \sqrt{q_j} \,\partial_{y}\phi(p_j,q_j)Z_{q_j} \Bigg) \\ \nonumber 
&=&  \frac{1}{m} \Bigg\{ \sum_{j\in D} q_j(1-q_j) 
\bigg[ \bigg( 1+\ln\frac{q_j}{p_j} \bigg)^2 
+ \Big(\frac{p_j}{q_j}\Big)^2 \bigg] \\ \nonumber 
&& - 2 \sum_{\{(i,j)\in D^2: \,i< j\}} 
q_i q_j \bigg[ \big(1+\ln\frac{q_i}{p_i} \big) 
\big(1+\ln\frac{q_j}{p_j} \big) 
+ \frac{p_i p_j}{q_i q_j} \bigg] \Bigg\}.
\end{eqnarray}
Using a slight modification of the arguments in \eqref{R_inequality} and \eqref{R_convergence}, it is verified the following convergence almost surely 
\begin{equation}\label{R_convergence_two}
    \Bigg|\sqrt{\frac{nm}{m\mathrm{V}_1(\mathbf{p},\mathbf{q})+n\mathrm{V}_2(\mathbf{p},\mathbf{q})}}\,\sum_{j \in D}R_j(n,m)\Bigg|\longrightarrow 0\quad \textrm{as}\quad n,m \longrightarrow \infty.
\end{equation}
Substituting equations \eqref{convergence_p}, \eqref{convergence_q} and \eqref{R_convergence_two} in \eqref{J_equality}, and since $Z_{\mathbf{p}}$ and $Z_{\mathbf{q}}$ are independent random vectors, and assuming that $\frac{nm}{n+m}\longrightarrow \gamma \in (0,1)$, finally we obtain 
\begin{equation}
    \sqrt{\frac{nm}{m\mathrm{V}_1(\mathbf{p},\mathbf{q})+n\mathrm{V}_2(\mathbf{p},\mathbf{q})}}\,\Big(J(\hat{\mathbf{p}}_n,\hat{\mathbf{q}}_m)-J(\mathbf{p},\mathbf{q})\Big)\overset{d}{\Longrightarrow}\mathcal{N}(0,1).
\end{equation}
\hfill$\square$
\vskip 0.2in

\section*{Declarations}
\textbf{Ethics approval and consent to participate. } This declaration is not applicable. \\ \\
\textbf{Consent for publication. } This declaration is not applicable. \\ \\
\textbf{Availability of data and materials. } The datasets supporting the findings of this study are available at the following repository: \href{https://github.com/Nicerova7/statistical_iv}{https://github.com/Nicerova7/statistical\_iv}. This repository provides access to the data used in the analysis, ensuring transparency and reproducibility of the results. Furthermore, we apply our test on fraud identification data and provide an open-source Python library, called ``statistical-iv" (\href{https://pypi.org/project/statistical-iv/}{https://pypi.org/project/statistical-iv/}). \\ \newline
\textbf{Competing Interests. }  The authors declare that they have no competing interests. \\ \\
\textbf{Funding. } No funding was received for conducting this study. \\ \\
\textbf{Authors' contributions. } All authors made significant contributions to the manuscript, including its design, data analysis, and drafting. All authors have read and approved the final version. \\ \\
\textbf{Acknowledgments. } This work has been partially supported by the Research Institute of the Faculty of Economics, Statistics and Social Sciences (IECOS) and the Vice-Rectorate of Research of the National University of Engineering (VRI-UNI) through the research project  25-2023-002666. \\ \\


\bibliography{references}

\end{document}